\documentclass[12pt]{article}

\usepackage{amsmath,amssymb,amsbsy,amsfonts,amsthm,latexsym,
                  amsopn,amstext,amsxtra,euscript,amscd}

\begin{document}

\newtheorem{theorem}{Theorem}
\newtheorem{lemma}[theorem]{Lemma}
\newtheorem{claim}[theorem]{Claim}
\newtheorem{cor}[theorem]{Corollary}
\newtheorem{prop}[theorem]{Proposition}
\newtheorem{definition}{Definition}
\newtheorem{question}[theorem]{Open Question}

\def\cA{{\mathcal A}}
\def\cB{{\mathcal B}}
\def\cC{{\mathcal C}}
\def\cD{{\mathcal D}}
\def\cE{{\mathcal E}}
\def\cF{{\mathcal F}}
\def\cG{{\mathcal G}}
\def\cH{{\mathcal H}}
\def\cI{{\mathcal I}}
\def\cJ{{\mathcal J}}
\def\cK{{\mathcal K}}
\def\cL{{\mathcal L}}
\def\cM{{\mathcal M}}
\def\cN{{\mathcal N}}
\def\cO{{\mathcal O}}
\def\cP{{\mathcal P}}
\def\cQ{{\mathcal Q}}
\def\cR{{\mathcal R}}
\def\cS{{\mathcal S}}
\def\cT{{\mathcal T}}
\def\cU{{\mathcal U}}
\def\cV{{\mathcal V}}
\def\cW{{\mathcal W}}
\def\cX{{\mathcal X}}
\def\cY{{\mathcal Y}}
\def\cZ{{\mathcal Z}}

\def\A{{\mathbb A}}
\def\B{{\mathbb B}}
\def\C{{\mathbb C}}
\def\D{{\mathbb D}}
\def\E{{\mathbb E}}
\def\F{{\mathbb F}}
\def\G{{\mathbb G}}
\def\I{{\mathbb I}}
\def\J{{\mathbb J}}
\def\K{{\mathbb K}}
\def\L{{\mathbb L}}
\def\M{{\mathbb M}}
\def\N{{\mathbb N}}
\def\O{{\mathbb O}}
\def\P{{\mathbb P}}
\def\Q{{\mathbb Q}}
\def\R{{\mathbb R}}
\def\S{{\mathbb S}}
\def\T{{\mathbb T}}
\def\U{{\mathbb U}}
\def\V{{\mathbb V}}
\def\W{{\mathbb W}}
\def\X{{\mathbb X}}
\def\Y{{\mathbb Y}}
\def\Z{{\mathbb Z}}

\def\ep{{\mathbf{e}}_p}

\def\scr{\scriptstyle}
\def\\{\cr}
\def\({\left(}
\def\){\right)}
\def\[{\left[}
\def\]{\right]}
\def\<{\langle}
\def\>{\rangle}
\def\fl#1{\left\lfloor#1\right\rfloor}
\def\rf#1{\left\lceil#1\right\rceil}
\def\le{\leqslant}
\def\ge{\geqslant}
\def\eps{\varepsilon}
\def\mand{\qquad\mbox{and}\qquad}

\def\vec#1{\mathbf{#1}}
\def\inv#1{\overline{#1}}
\def\vol#1{\mathrm{vol}\,{#1}}

\newcommand{\comm}[1]{\marginpar{%
\vskip-\baselineskip 
\raggedright\footnotesize
\itshape\hrule\smallskip#1\par\smallskip\hrule}}

\def\xxx{\vskip5pt\hrule\vskip5pt}


\title{\bf Approximation  by Several Rationals}

\author{ 
{\sc Igor E. Shparlinski} \\
{Department of Computing, Macquarie University} \\
{Sydney, NSW 2109, Australia} \\
{igor@ics.mq.edu.au}}

\date{\today}
\pagenumbering{arabic}

\maketitle

\begin{abstract} Following T.~H.~Chan, we consider 
the problem of approximation  of a given rational
fraction $a/q$ by sums of several rational fractions $a_1/q_1, \ldots,
 a_n/q_n$ with smaller denominators. We show that in the special
cases of $n=3$  and $n=4$ and certain admissible ranges for the 
denominators $q_1,\ldots, q_n$,  one can  improve  a
result of T.~H.~Chan by using a different approach. 
\end{abstract}

 \paragraph*{2000 Mathematics Subject Classification:} 	  11J04, 11N25

\section{Introduction}
\label{sec:intro}

T.~H.~Chan~\cite{Chan} has recently considered the question 
of approximating real numbers by sums of several 
 rational fractions $a_1/q_1, \ldots,
 a_n/q_n$ with bounded  denominators. 

In the  special case
of  $n = 3$ the result of T.~H.~Chan~\cite{Chan} can be reformulated 
as follows. Given two integers $a$ and $q\ge 1$, 
for any   $Q \ge q$ 
there are integers  $a_i$ and $q_i$ with 
 $1 \le q_i \le Q^{4/7+o(1)}$, $i =  1, 2, 3$, and  
such that  
$$
\left| \frac{a}{q} - \frac{a_1}{q_1} - \frac{a_2}{q_2} - 
\frac{a_3}{q_3}\right|\le \frac{1}{qQ^{1+o(1)}}.
$$
We remark that the numerators  $a_1, a_2, a_3$ 
can be negative.

In this  paper   we use   different  approach to show that when 
$Q$ is large enough, namely, when $Q\ge q^{2+\varepsilon}$ the same result holds with
$1/3$ instead of $4/7$. We also obtain more explicit constants. 

Similarly, for $n=4$, we see from~\cite{Chan} 
that for any    $Q \ge q$ 
 there are integers  $a_i$ and $q_i$ with 
 $1 \le q_i \le Q^{2/5+o(1)}$, $i =  1, 2, 3,4$, and  
such that  
$$
\left| \frac{a}{q} - \frac{a_1}{q_1} - \frac{a_2}{q_2} - 
\frac{a_3}{q_3}- \frac{a_4}{q_4}\right|\le \frac{1}{qQ^{1+o(1)}}.
$$
In this case, under the same condition  $Q\ge q^{2+\varepsilon}$ 
we replace 
$2/5$ instead of $1/4$.

Our approach is based on a result of~\cite{Shp2} about the uniformity 
of distribution in residue of rather general products. More precisely,
it is shown in~\cite{Shp2} that for any set $\cX \in [1,X]$ of integers
$x$ with $\gcd(x,q)=1$ and
for any interval $[Z+1, Z+Y]$, for the number $M_{u,q}(\cX;Y,Z)$ of solutions
to  the congruence 
$$
u \equiv xy \pmod q, \qquad x \in \cX, \ y \in [Z+1, Z+Y], 
$$
we have  
\begin{equation}
\label{eq:Aver Prod}
\sum_{u=1}^q  \left|M_{u,q}(\cX;Y,Z) -   \# \cX  \frac{Y}{q}  \right|^2 
\le \# \cX (X+Y)  q^{o(1)}.
\end{equation}


\section{Approximation  by Three Rationals}

\begin{theorem}
\label{thm:Approx 3}  Let $a$ and $q\ge 1$ be integers with 
$\gcd(a, q) = 1$.    For any fixed $\varepsilon >0$ and sufficiently large $q$, 
for any integer  $Q \ge q^{2+\varepsilon}$ 
 there are integers  $a_i$ and $q_i$ with 
$1 \le q_i \le 2 Q^{1/3}$, $i =  1, 2, 3$, and  
such that
$$
\left| \frac{a}{q} - \frac{a_1}{q_1} - \frac{a_2}{q_2} - 
\frac{a_3}{q_3}\right|\le \frac{1}{qQ } 
$$ 
holds. 
\end{theorem}

\begin{proof} We note that it is enough to
show that there are positive integers
$ q_1, q_2, q_3 \le  2Q^{1/3}$ with 
\begin{equation}
\label{eq:Prod 3}
q_1  q_2  q_3 \ge Q
\end{equation}
such that
\begin{equation}
\label{eq:gcd 3}
\gcd( q_1, q_2) = \gcd( q_1,  q_3) = \gcd( q_2,  q_3) = 1   
\end{equation}
and
\begin{equation}
\label{eq:cong}
aq_1  q_2  q_3 \equiv 1   \pmod q.
\end{equation}

Indeed, from~\eqref{eq:cong} we conclude that
$aq_1  q_2  q_3 = 1 + bq$
for some integer  $b$. Since~\eqref{eq:gcd 3}  implies that
$$
\gcd( q_1  q_2,  q_1  q_3, q_2 q_3) = 1  
$$
then
$$
b = a_1q_2q_3 + a_2q_1q_3 + a_3q_1q_2
$$
for are some integers   $a_1,a_2,a_3$.  Thus 
$$
\left| \frac{a}{q} - \frac{a_1}{q_1} - \frac{a_2}{q_2} - 
\frac{a_3}{q_3}\right| = \frac{1}{qq_1q_2q_3}
\le \frac{1}{qQ}.
$$
Let us put  $R = \fl{2 Q^{1/3}}$. We  may assume that $R < q$
since otherwise we simply choose $a_1=1$, $a_2 = a_3 = 0$, 
$q_1= q$, $q_2 = q_3 = 1$.

 We now consider
\begin{itemize}
\item the set $\cS$ consisting of integers  $s \in [R/3, R/2)$;
\item the set $\cP$ consisting of primes $p \in [R/2, 3R/4)$ with $\gcd( p, q) = 1$;
\item the set $\cL$ consisting of primes  $\ell \in [3R/4, R]$
 with $\gcd( \ell, q) = 1$.
\end{itemize}

Since $q$ may have at most $O(\log q)$ prime divisors, by the prime 
number theorem we see that 
$$
\# \cS, \#\cP, \#\cL \ge R^{1 + o(1)}.
$$
Clearly, if we  take $q_1 = s \in \cS$, $q_2 = p \in \cP$ and $q_3 = \ell \in \cL$ 
then~\eqref{eq:gcd 3}  is satisfied and we also have~\eqref{eq:Prod 3}.   Thus  it
is enough to   show that the congruence
$$
sp\ell \equiv 1 \pmod q, \qquad s \in \cS,\ p \in \cP,\ \ell \in \cL,
$$ 
has a solution.
For an integer  $u \in [1 , q]$ we denote
by  $N(u)$ the number   of solutions to the
congruence 
\begin{equation}
\label{eq:cong u}
sp  \equiv  u \pmod q, \qquad s \in \cS,\  p \in \cP.
\end{equation} 
Let $\cU$ be the set  of integers $u \in [1, q]$ for which 
the above congruence has a solution, that is, $N(u) > 0$.
It is enough to show that the congruence
\begin{equation}
\label{eq:target}
u\ell \equiv 1 \pmod q, \qquad u \in\cU, \ \ell \in \cL,
\end{equation}
has a solution.

Also let $\cV$ be the set of remaining integers $u \in [1 , q]$ 
with $N(u)= 0$. It follows from~\cite{Shp2}, see~\eqref{eq:Aver Prod}, that
$$
\sum_{u=1}^q \left|N(u) - \frac{\#\cS \#\cP}{q}\right|^2 \le R^2 q^{o(1)}. 
$$
Hence
$$
\# \cV \(\frac{\#\cS \#\cP}{q}\)^2 \le R^2 q^{o(1)}
$$
which implies that $\# \cV  \le  R^{-2} q^{2 + o(1)}$. 
Recalling that 
 $R \ge  2 Q^{1/3} -1 \ge  q^{2/3 + \varepsilon/3}$, we see that 
$$
 \# \cL - \# \cV = R^{1 +o(1)} - R^{-2} q^{2 + o(1)} > 0
$$ 
provided that $q$ is large enough. Therefore
 the congruence~\eqref{eq:target} has  a
solution, which concludes the proof. 
 \end{proof} 

\section{Approximation  by Four Rationals}

We now use a similar approach for  approximations by four rational 
fractions. 

\begin{theorem}
\label{thm:Approx 4}  Let $a$ and $q\ge 1$ be integers with 
$\gcd(a, q) = 1$.  For any fixed $\varepsilon >0$ 
and sufficiently large $q$,
for any integer  $Q \ge q^{2+\varepsilon}$   
 there are integers  $a_i$ and $q_i$ with 
 $1 \le q_i \le 2 Q^{1/4}$, $i =  1, 2, 3$, and  
such that
$$
\left| \frac{a}{q} - \frac{a_1}{q_1} - \frac{a_2}{q_2} - 
\frac{a_3}{q_3}- \frac{a_4}{q_4}\right|\le \frac{1}{qQ } 
$$ 
holds. 
\end{theorem}

\begin{proof} We proceed as in the proof of Theorem~\ref{thm:Approx 3}. 
In particular, we see that   it
is enough to show that there are positive integers
$ q_1, q_2, q_3, q_4\le  2Q^{1/4}$ with 
\begin{equation}
\label{eq:Prod 4}
q_1  q_2  q_3 q_4 \ge Q
\end{equation}
such that
\begin{equation}
\label{eq:gcd 4}
\gcd( q_i, q_j) = 1 , \qquad 1 \le i < j\le 4, 
\end{equation}
and
$$
aq_1  q_2  q_3 q_4 \equiv 1   \pmod q.
$$

Let us put  $R = \fl{2 Q^{1/4}}$. As before, we remark 
that we  may assume that $R < q$ 
since otherwise the result is trivial. 

 We now consider
\begin{itemize}
\item the set $\cS$ consisting of integers  $s \in [R/4, R/3)$;
\item the set $\cP$ consisting of primes $p \in [R/3,  2R/3)$ with $\gcd( p, q) = 1$;
\item the set $\cL$ consisting of primes  $\ell \in [2R/3, 3R/4)$
 with $\gcd( \ell, q) = 1$;
\item the set $\cR$ consisting of primes  $r \in [3R/4, R]$
 with $\gcd( r, q) = 1$.
\end{itemize}
Again,  by the prime number theorem we have:
$$
\# \cS, \#\cP, \#\cL, \#\cR \ge R^{1 + o(1)}.
$$
Clearly, if we  take $q_1 = s \in \cS$, $q_2 = p \in \cP$, $q_3 = \ell \in \cL$ 
and $q_4 = r \in \cR$ 
then~\eqref{eq:gcd 4}  is satisfied and we also have~\eqref{eq:Prod 4}.   Thus  it
is enough to   show that the congruence
$$
sp\ell r \equiv 1 \pmod q, \qquad s \in \cS,\ p \in \cP,\ \ell \in \cL, \ r \in \cR,
$$ 
has a solution.

As in the proof of Theorem~\ref{thm:Approx 3} we note the
set $\cV$ of integers $u \in [1 , q]$ for which the
congruence~\eqref{eq:cong u} does not have a solution is 
of cardinality $\# \cV  \le  R^{-2} q^{2 + o(1)}$.

Let $\cW$ be the set  of integers $w \in [1 , q]$ 
which are of the form
$w  \equiv \ell r \pmod q$ with $\ell \in \cL$ and $r \in \cR$. 
We note that 
$\#\cL \#\cR = R^{2+o(1)}$  
products $\ell r$ are distinct integers in the interval $[1, R^2]$. 
Since there are at most $R^2/q +1$ integers $t \in [1, R^2]$
in the same residue class modulo $q$, 
we obtain 
$$
\# \cW \ge R^{2+o(1)}\(R^2/q +1\)^{-1}.
$$
Since $R^2 \ge(2 Q^{1/4}-1)^2 \ge Q^{1/2} \ge q^{1+\varepsilon/2}$
(provided $q$ is large enough) we see that $R^2/q +1 \le 2R^2/1$.
Hence $\# \cW = q^{1 + o(1)}$. 
We now see that 
$$\#\cW - \# \cV = q^{1 + o(1)} - R^{-2} q^{2 + o(1)} 
> 0
$$ 
provided that $q$ is large enough. The desired 
result now follows. 
\end{proof} 

We remark that in both Theorem~\ref{thm:Approx 3}
and Theorem~\ref{thm:Approx 4} the coefficient $2$ 
in the bound on the denominators can be replaced
by  any constant $c > 1$. 

%

\section{Comments}

It is  natural to try to use~\eqref{eq:Aver Prod} 
to improve  the corresponding bound from~\cite{Chan}
for larger values of $n$ too. Although some results can 
be obtained on this way, for $n \ge 5$  we have not been able
to achieve this. In fact, it seems quite plausible that for $n \ge 5$,
instead of using the bound~\eqref{eq:Aver Prod}  from~\cite{Shp2}, 
one can study the solvability of the congruence
$$
q_1\ldots q_n \equiv 1 \pmod q
$$
with ``small'' $q_1,\ldots,q_n$ by using bounds of character sums
in the same style as in~\cite{Shp1,Shp3}.

\section*{Acknowledgement}

The author is grateful to  Tsz Ho Chan for useful discussions.

This work was supported in part by ARC grant
DP0556431.

\end{document}